\documentclass{amsart}
\usepackage{graphics}
\usepackage{epsfig}

\title{Catalan Numbers and Schubert Polynomials for $w=1(n+1)\cdots2$}
\author{Alexander Woo}
\thanks{Research supported in part by NSF grant DMS-0301072.}

\newenvironment{smallrc}{\begin{trivlist}\item\centering\small$}
			{$\end{trivlist}}

\newenvironment{tinygraph}{\begin{trivlist}\item\centering\tiny$}
			{$\end{trivlist}}

\def\textcross{\ \smash{\lower2pt\hbox{\rlap{\hskip2.07pt\vrule height10pt}}
                \raise2.8pt\hbox{\rlap{\hskip-3pt \vrule height.4pt depth0pt
                width10.7pt}}}\hskip10.7pt\!\!}
\def\textelbow{\ \hskip.1pt\smash{\raise2.8pt%
                \hbox{\co \hskip 2.07pt\rlap{\rlap{\char'005} \char'007}
                \lower5.2pt\rlap{\vrule height1.5pt}
                \raise4.2pt\rlap{\vrule height1.5pt}}
                \raise2.8pt\hbox{%
                  \rlap{\hskip-5.25pt \vrule height.4pt depth0pt width1.5pt}%
                  \rlap{\hskip4.25pt \vrule height.4pt depth0pt width1.5pt}}}
                \hskip6.7pt}

\font\co=lcircle10
\def\petit#1{{\scriptstyle #1}}

\def\jr{\smash{	\raise2pt\hbox{\co \rlap{\rlap{\char'005} \char'007}}
		\raise6pt\hbox{\rlap{\vrule height2pt}}
		\raise2pt\hbox{\rlap{\hskip4pt \vrule height0.4pt depth0pt
                 width2.5pt}}
		\raise2pt\hbox{\rlap{\hskip-6pt \vrule height.4pt depth0pt
                 width2.2pt}}
		\lower4pt\hbox{\rlap{\vrule height2.5pt}}}}
\def\je{\smash{\raise2pt\hbox{\co \rlap{\rlap{\char'005}
                \phantom{\char'007}}}\raise6pt\hbox{\rlap{\vrule height2pt}}
	       \raise2pt\hbox{\rlap{\hskip-6pt \vrule height.4pt depth0pt
                width2.2pt}}}}
\def\+{\smash{\lower4pt\hbox{\rlap{\vrule height12pt}}
                \raise2pt\hbox{\rlap{\hskip-6pt \vrule height.4pt depth0pt
                width12.5pt}}}}
\def\er{\smash{	\raise2pt\hbox{\co \rlap{\phantom{\rlap{\char'005}} \char'007}}
		\raise6pt\hbox{\rlap{\phantom{\vrule height2pt}}}
		\raise2pt\hbox{\rlap{\hskip4pt \vrule height0.4pt depth0pt
                 width2.5pt}}
		\raise2pt\hbox{\rlap{\phantom{%
		 \hskip-6pt \vrule height.4pt depth0pt width2.2pt}}}
		\lower4pt\hbox{\rlap{\vrule height2.5pt}}}}
\def\hor{\smash{\lower2pt\hbox{\rlap{\phantom{\vrule height10pt}}}
                \raise2pt\hbox{\rlap{\hskip-6pt \vrule height.4pt depth0pt
                width12.5pt}}}}
\def\ver{\smash{\lower2pt\hbox{\rlap{\vrule height10pt}}
                \raise2pt\hbox{\rlap{\phantom{%
		\hskip-6pt \vrule height.4pt depth0pt width12.5pt}}}}}

\def\perm#1#2{\hbox{\rlap{$\petit {#1}_{\scriptscriptstyle #2}$}}%
                \phantom{\petit 1}}

	{%
	 \def\*{\makebox[0ex]{\footnotesize$+\,$}}%
	 \begin{array}{|*{#1}{@{\ \;\:}c|}}}
	{\end{array}}

\newenvironment{pipedream}[1]%
	{%
	 \begin{array}{c*{#1}{@{\ \ \;}c}}}
	{\end{array}}

\newenvironment{shortstretchpipedream}[1]%
	{%
	 \begin{array}{c*{#1}{@{\ \ \;}c}}}
	{\end{array}}

\def\adots{{.\hspace{0.5pt}\raisebox{3pt}{.}\hspace{0.5pt}\raisebox{6pt}{.}}}

\def\shortvdots{\smash\vdots}

\newtheorem{prop}{Proposition}

\begin{document}

\begin{abstract}
We show that the Schubert polynomial $\mathfrak{S}_w$ specializes to
the Catalan number $C_n$ when $w=1(n+1)\cdots2$.  Several proofs of
this result as well as a $q$-analog are given.  An application to the
singularities of Schubert varieties is given.
\end{abstract}

\maketitle

\section{Introduction}

The Catalan numbers $C_n=\frac{1}{n+1}\binom{2n}{n}$ are ubiquitous in
combinatorics.  Among other things they count Dyck paths, which are
lattice paths from $(0,0)$ to $(n,n)$ staying above the main diagonal,
and rooted binary trees with $n$ leaves (where each internal node has
exactly 2 children).  In this paper we present a perviously
undiscovered connection between Catalan numbers and certain Schubert
polynomials.  Our result also has an interpretation in terms of the
geometry of certain Schubert varieties.

Let $w\in S_n$ be a permutation.  The Schubert polynomial
$\mathfrak{S}_w$ is a polynomial representative for the class of the
Schubert variety $X_{w_0w}$ in the flag manifold $GL(n)/B$.  A
combinatorial formula for Schubert polynomials, which for our purposes
can be taken as their definition, was proven independently in
\cite{BJS} and \cite{FS}.  Diagrams illustrating this formula, called
rc-graphs, were introduced in \cite{FK} and extensively studied in
\cite{BB}.

For the remainder of this paper, let $w_n$ denote the permutation
$1(n+1)n\cdots2$ in $S_{n+1}$.  Using rc-graphs, we will show that
that the principal specialization of the Schubert polynomial for this
permutation, $\mathfrak{S}_{w_n}(1,q,\cdots,q^n)$, is equal to
$q^{\binom{n}{3}}C_n(q)$, where $C_n(q)$ is the Carlitz-Riordan
$q$-analogue of the Catalan numbers originally introduced in
\cite{CR}.  Section \ref{sect:rc-graphs} gives the definition of
rc-graphs and proves our main theorem via a recurrence counting
rc-graphs for $w_n$.  Section \ref{sect:elem-bij} gives a bijection
to Dyck paths, and Section \ref{sect:edel-green} shows that this
bijection can also be described using the Edelman-Greene
correspondence.  Since $w_n$ is its own inverse in $S_{n+1}$,
transposition is a natural involution on its rc-graphs.  In Section
\ref{sect:tree-inv}, we describe our bijection in terms of binary
trees, giving a correspondence between features of the rc-graph known
as elbow joints and the internal nodes of a binary tree.  From this
description it will be evident that transposing an rc-graph
corresponds to flipping a binary tree around its vertical axis.

Finally we discuss in Section \ref{sect:geom} the geometric example
which originally motivated this study.  Let $w^\prime_n$ be the
permutation $(n+2)23\cdots(n+1)1$ in $S_{n+2}$.  We show that the
multiplicity of the Schubert variety $X_{w^\prime_n}$ at its most
singular point is given by $C_n$.  Experimental evidence \cite{me}
suggests that this is higher than the multiplicity of any point on any
other Schubert variety of $GL(n+2)/B$.

\section{RC-Graphs}\label{sect:rc-graphs}

\begin{figure}
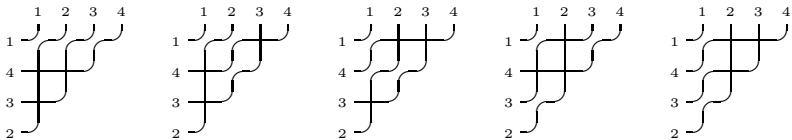


\begin{tinygraph}
\begin{pipedream}{4}
       &\perm1{}&\perm2{}&\perm3{}&\perm4{}\\
\petit1&   \jr  &   \jr  &   \jr  &  \je   \\
\petit4&   \+   &   \+   &   \je  &\\
\petit3&   \+   &   \je  &        &\\
\petit2&   \je  &        &        &\\
\end{pipedream}
\quad
\begin{pipedream}{4}
       &\perm1{}&\perm2{}&\perm3{}&\perm4{}\\
\petit1&   \jr  &   \jr  &   \+    &  \je   \\
\petit4&   \+   &   \jr  &   \je  &\\
\petit3&   \+   &   \je  &        &\\
\petit2&   \je  &        &        &\\
\end{pipedream}
\quad
\begin{pipedream}{4}
       &\perm1{}&\perm2{}&\perm3{}&\perm4{}\\
\petit1&   \jr  &   \+   &   \+    &  \je   \\
\petit4&   \jr &    \jr  &   \je  &\\
\petit3&   \+   &   \je  &        &\\
\petit2&   \je  &        &        &\\
\end{pipedream}
\quad
\begin{pipedream}{4}
       &\perm1{}&\perm2{}&\perm3{}&\perm4{}\\
\petit1&   \jr  &   \+   &   \jr  &  \je   \\
\petit4&   \+   &   \+   &   \je  &\\
\petit3&   \jr  &   \je  &        &\\
\petit2&   \je  &        &        &\\
\end{pipedream}
\quad
\begin{pipedream}{4}
       &\perm1{}&\perm2{}&\perm3{}&\perm4{}\\
\petit1&   \jr  &   \+   &   \+   &  \je   \\
\petit4&   \jr  &   \+   &   \je  &\\
\petit3&   \jr  &   \je  &        &\\
\petit2&   \je  &        &        &\\
\end{pipedream}
\end{tinygraph}

\caption{The rc-graphs for $w = 1432$}\label{fig:rc-graphs}
\end{figure}

Let $w\in S_n$ be a permutation.  An {\it rc-graph} for $w$ is a
filling of the upper-left half of an $n \times n$ array with cross
pieces ($\textcross$) and elbow joints ($\textelbow$) such that the
strand entering the left in row $i$ exits the top in column $w(i)$,
with the additional condition that no two strands cross more than
once.  This second condition can alternatively be stated as there
being exactly $l(w)$ cross pieces.  For example, the 5 rc-graphs for
$w=1432$ are shown in figure \ref{fig:rc-graphs}.  It was shown
independently in \cite{BJS} and \cite{FS} that rc-graphs are related
to Schubert polynomials by the formula
$$\mathfrak{S}_w(x_1,\cdots, x_n) = \sum_{D\in{\mathcal{RC}(w)}}
\prod_{(i,j)\in C(D)} x_i,$$ where $\mathcal{RC}(w)$ is the set of
rc-graphs for $w$, and $C(D)$ are the locations of the cross pieces
in $D$ (indexed so that $(1,3)$ would be in $C(D)$ if $D$ has a cross
in the top row and the third column).  For example, we have
$\mathfrak{S}_{1432}=x_2^2x_3 + x_1x_2x_3 + x_1^2x_3 + x_1x_2^2 +
x_1^2x_2$, with the 5 terms corresponding to the 5 rc-graphs in the
figure \ref{fig:rc-graphs} from left to right.

First we will show that there are in fact Catalan many rc-graphs for
$w_n$; this will allow us to prove that the combinatorial maps we give
from the set of rc-graphs for $w_n$ to other Catalan objects are in
fact bijections by only showing that they are injections or
surjections.

The Carlitz-Riordan $q$-Catalan numbers $C_n(q)$ are defined by the
recurrence $C_n(q)= \sum_{k=0}^{n-1} q^kC_{n-k-1}(q)C_k(q)$, with
$C_0(q)=1$; under the interpretation of Catalan numbers as counting
partitions $\lambda$ whose Young diagrams fit inside the Young diagram
of the staircase partition $\delta_n=n-1,\cdots,1$, we have
$C_n(q)=\sum_\lambda q^{\binom{n}{2}-\left|\lambda\right|}$.

\begin{prop}

$$ \mathfrak{S}_{w_n}(1,q,\cdots,q^n) = q^{\binom{n}{3}}C_n(q).$$
\label{prop:count}
\end{prop}

\begin{proof}

We proceed by induction on $n$.  The proposition is clear for $n=1$.

Given a rc-graph $D$, the strand $D_{n+1}$ entering the left in row
$n+1$ and exiting the top in column $w_n(n+1)=2$ travels through only
the first and second columns.  There is a unique $k$, $1\leq k\leq n$,
such that $D_{n+1}$ goes through both the square $(k,1)$ and $(k,2)$.

Fixing $k$, there must be cross pieces at $(i,j)$ for all $i,j$ with
$1\leq i\leq k-1$, $2\leq j\leq n+2-k$; this is true by definition for
$j=2$, and each strand $D_l$ for $k+1\leq l\leq n$ crosses $D_{n+1}$
in row $l$ and never travels to the right of column $w(l)=n+3-l$, so they
must each go straight through the topmost $k-1$ places in column
$n+3-l$, which therefore must all be cross pieces.  Note each column
of these crosses contributes $\binom{k-1}{2}$ to the weight, for a
total of $(n+1-k)\binom{k-1}{2}$.

We must also have cross pieces at $(i,1)$ for $k+1\leq j\leq n$, and
these contribute $\binom{n}{2}-\binom{k}{2}$ to the weight.

Now let $D^\prime$ be the portion of $D$ consisting of squares $(i,j)$
with $k\leq i\leq n$ and $2\leq j\leq n+2-k$; this is an rc-graph for
$w_{n-k}$.  Also, let $D^{\prime\prime}$ be the portion of $D$
consisting of the squares $(i,j)$, where $1\leq i\leq k$ and either
$j=1$ or $n+3-k\leq j\leq n+1$, ignoring the intervening cross pieces
in columns $2$ through $n+2-k$; this is an rc-graph for $w_{k-1}$.

\begin{figure}
  \epsfig{file=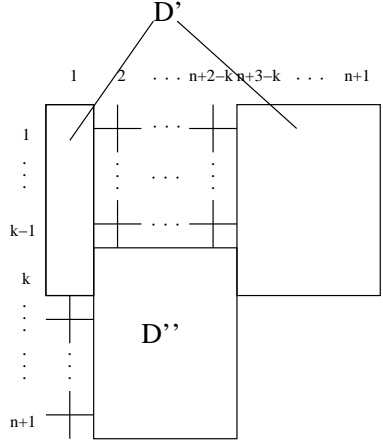}
\caption{Breaking up the rc-graph $D$}\label{fig:rc-split}
\end{figure}

Figure \ref{fig:rc-split} illustrates the situation; the elbow joints
have been left out of the diagram to make it smaller and more
readable.  The numbers on the left denote in this figure the row
number; note that $w_n(k)=n+3-k$.

For convenience, let $\operatorname{wt}(D)=\sum_{(i,j)\in C(D)}
(i-1)$.  Note $\operatorname{wt}(D)$ is the exponent of $q$ in the
term of $\mathfrak{S}_w(1,q,\cdots, q^n)$ corresponding to the
rc-graph $D$.  The rc-graph $D^{\prime\prime}$ contributes exactly
$\operatorname{wt}(D^{\prime\prime})$ to $\operatorname{wt}(D)$.  The
cross pieces in $D^\prime$, however, are each shifted $k-1$ rows down,
and there are $\binom{n-k}{2}$ cross pieces in $D^\prime$, so
$D^\prime$ contributes
$\operatorname{wt}(D^\prime)+(k-1)\binom{n-k}{2}$ to the weight of
$D$.

Therefore, for an rc-graph $D$ for $w_n$, we have
$$wt(D) = wt(D^{\prime\prime})
        + (k-1)\binom{n-k}{2}+\operatorname{wt}(D^\prime)
        + (n+1-k)\binom{k-1}{2}
        + \binom{n}{2}-\binom{k}{2}.$$

Let $F_n(q)$ denote $\mathfrak{S}_{w_n}(1,q,\cdots,q^n) =
\sum_{D\in\mathcal{RC}(w_n)} q^{\operatorname{wt}(D)}$.

Then we have
\begin{align}
F_n(q) &= \sum_{k=1}^{n} q^{(k-1)\binom{n-k}{2}+(n-k+1)\binom{k-1}{2}+\binom{n}{2}-\binom{k}{2}}F_{k-1}(q)F_{n-k}(q) \notag\\
         &= \sum_{k=1}^{n} q^{(k-1)\binom{n-k}{2}+(n-k-1)\binom{k-1}{2}+\binom{n}{2}-\binom{k}{2}+\binom{k-1}{3}+\binom{n-k}{3}}C_{k-1}(q)C_{n-k}(q) \notag\\
	 &= \sum_{k=1}^{n} q^{\binom{n}{3}+n-k}C_{k-1}(q)C_{n-k}(q) \notag\\
         &= q^{\binom{n}{3}} \sum_{k=0}^{n-1} q^kC_{n-k-1}(q)C_k(q) \notag\\
	 &= q^{\binom{n}{3}} C_n(q), \notag
\end{align}
which proves the proposition.

\end{proof}

\section{A Bijection to Dyck Paths}\label{sect:elem-bij}

Given proposition \ref{prop:count}, we would like bijections between
rc-graphs for $w_n$ and other sets of objects counted by Catalan
numbers.  One such set is the set of partitions $\lambda$ whose Young
diagrams fit inside the Young diagram for the staircase partition
$\delta_n=n-1,\cdots,1$, or, equivalently, such that $\lambda_k\leq
n-k$ for all $k$.  We denote this set $\mathcal{DP}(n)$.

Let $D$ be an rc-graph for $w_n$, and let $E(D)$ be the set of
locations of its elbow joints.  We can then associate a partition
$\lambda(D)$ to $D$ by requiring that the parts of its conjugate
$\lambda^\prime(D)$ be, as a multiset, $\{\{j-1|(i,j)\in E(D),
(i,j)\neq(0,0)\}\}$; in other words, each elbow joint in the $i$-th
row of $D$ should contribute a part of $i-1$ to $\lambda^\prime(D)$.

The proof will involve generalized inverse chute moves, which are
local moves first given in \cite{BB} that, given an rc-graph for some
permutation $w$, allows one to generate new rc-graphs for $w$.  Let
$D$ be any rc-graph for some permutation $w$.  Suppose $D$ has an
elbow joint at $(i,j)$, and that the following all hold for some
$i^\prime>i$ and some $j^\prime<j$:
\begin{enumerate}
\item For each $k$, $i<k<i^\prime$, $(k,j)$ is a cross piece, and $(i^\prime,j)$ is an elbow joint.
\item For each $k$, $j^\prime<k<j$, $(i,k)$ is a cross piece, and $(i,j^\prime)$ is an elbow joint.
\item For each $k$, $i<k\leq i^\prime$, $(k,j^\prime)$ is a cross piece.
\item For each $k$, $j^\prime\leq k<j$, $(i^\prime,k)$ is a cross piece.
\end{enumerate}
Then the diagram $D^\prime$ obtained from $D$ by changing the elbow
joint at $(i,j)$ to a cross piece and the cross piece at
$(i^\prime,j^\prime)$ to an elbow joint is also an rc-graph for $w$,
since in both cases the strand entering this area from the left at
$(i^\prime,j^\prime)$ exits $(i,j)$ on the right, the strand entering
from the bottom at $(i^\prime,j^\prime)$ exits $(i,j)$ on top, all
other strands are unchanged, and the number of cross pieces remains
the same.

\begin{prop}
\label{prop:elem-bij}
The above described map from $D$ to $\lambda(D)$ gives a bijection
from $\mathcal{RC}(w_n)$ to $\mathcal{DP}(n)$.
\end{prop}

\begin{proof}
We will show that this map is surjective and then appeal to
proposition \ref{prop:count}.

\begin{figure}
\begin{smallrc}
\begin{shortstretchpipedream}{6}
             & \perm 1{} & \perm 2{} & \perm\cdots{} & \perm\cdots{} & \perm {n}{} & \perm {n+1}{}\\
\petit 1     &    \jr    &    \jr    & \petit\cdots &  \jr   &  \jr        & \je \\
\petit {n+1} &    \+     &    \+     & \petit\cdots &   \+   &  \je    &      \\
\petit {n}   &    \petit\shortvdots &   \petit\shortvdots  & \petit\adots &  \je   &      & \\
\petit \shortvdots       &    \petit\shortvdots &   \+      & \petit\adots &        &      & \\
\petit \shortvdots       &    \+     &   \je     &        &        &      & \\
\petit 2     &    \je    &           &        &        &      & \\
\end{shortstretchpipedream}
\end{smallrc}

\caption{The pipe dream $D_\mathrm{bot}$}\label{fig:rc-bottom}
\end{figure}

Begin with the bottom rc-graph for $w_n$ with elbow joints in the
first row and cross pieces everywhere else as shown in figure
\ref{fig:rc-bottom}; we denote this rc-graph $D_\mathrm{bot}$.  Note
that $\lambda(D_\mathrm{bot})=\emptyset$.  Now we construct the
rc-graph $D$ with $\lambda(D)=\lambda$.  Let $\lambda^\prime_i$ denote
the $i$-th part of $\lambda^\prime$, with the parts in decreasing
order.  Now, for each $i$ starting from 1, if $\lambda^\prime_i=k$,
take the rightmost cross piece in row $k+1$ that is not under a cross
piece in the top row, and let $l$ be the column it is in.  There
exists such a cross piece since, as $\lambda_j<n-j$, $\lambda^\prime$
has at most $n-j$ parts of size greater than or equal to $j$, and
$D_\mathrm{bot}$ has $n-j$ cross pieces in row $j+1$.  Now let $m$ be
the leftmost column to the right of $l$ such that $(1,l)$ has an elbow
joint.  Then moving the cross piece at $(k+1,l)$ to $(1,m)$ is a
generalized chute move as follows.  Condition 2 follows from the
definition of $m$, and conditions 1 and 3 follow from $D_\mathrm{bot}$
having only cross pieces in rows 2 through $k$, none of which have
been moved.  Condition 4 holds because the leftmost elbow joint in row
$k+1$ to the right of column $l$ must either be at the end of the row,
or the result of an immediately preceding move from row $k+1$, and
therefore in column $m$.  Performing these corresponding generalized
chute moves for all the parts of $\lambda^\prime$ in decreasing order
constructs an rc-graph for $w_n$ that goes to $\lambda$ under the
given map.  Therefore the map is a surjection and hence by our earlier
count of rc-graphs a bijection.

\end{proof}

\section{Edelman-Greene Correspondence}\label{sect:edel-green}

Our bijection has a second description in terms of the Edelman-Greene
correspondence, a variant of the usual RSK correspondence, originally
introduced in \cite{EG} and extended to the semi-standard case used
here in \cite{BJS}. This correspondence associates to each rc-graph a
pair $(P, Q)$ of column-strict Young tableaux of the same shape in
such a way that if two rc-graphs have the same $P$-tableau, they must
be rc-graphs for the same permutation.

The Edelman-Greene correspondence works as follows.  First convert the
rc-graph into a sequence $((a_1,\alpha_1), \cdots, (a_{l(w)},
\alpha_{l(w)}))$ of pairs of numbers to put into the tableaux as
follows.  Reading each row of the rc-graph from left to right and
starting with the top row, if the $k$-th cross piece is encountered at
$(i_k,j_k)$, let $a_k=i_k$ and $\alpha_k=i_k+j_k$.  Then we insert the
$\alpha_k$ one by one to create the $P$ tableau using Edelman-Greene
insertion; this is identical to RSK insertion except that, when
inserting the letter $i$ into a row with both an $i$ and an $i+1$,
that row remains unchanged and an $i+1$ is bumped into the next row.
As in RSK insertion, after all the bumping associated with inserting a
single letter $\alpha_k$ is completed, $a_k$ is added to the $Q$
tableau so that it has the same shape as $P$.  When we are finished,
$P$ will be both row and column-strict, but $Q$ will only be
row-strict.  Since column-strict tableaux are customary, we transpose
both tableaux.

\begin{figure}

$\begin{array}{ccccc}
1 & 2 & \cdots & \cdots & n \\
2 & \ddots & \ddots & n & \\
\vdots & \ddots & n & & \\
\vdots & n &  &  & \\
n &  &  &  &  \\
\end{array}$

\caption{$P$ tableau for $w=n\cdots1$}\label{fig:longwordEG}
\end{figure}

\begin{figure}

$\begin{array}{ccccc}
2 & 3 & \cdots & \cdots & n+1 \\
3 & \ddots & \ddots & n+1 & \\
\vdots & \ddots & n+1 & & \\
\vdots & n+1 &  &  & \\
n+1 &  &  &  &  \\
\end{array}$

\caption{$P$ tableau for $w_n$}\label{fig:ourwordEG}
\end{figure}

Edelman and Greene showed that the $P$ tableau for the long word
$n\cdots1$ is always the one in figure \ref{fig:longwordEG}.  As a trivial
corollary, the $P$ tableau for any rc-graph for $w_n$ must be the one
in figure \ref{fig:ourwordEG}.  In the case of the long word, they also
gave an inverse to their insertion procedure similar to {\it
evacuation}, an operation on tableaux originally due to
Sch\"utzenberger\cite{S}.  Take the $Q$-tableau, and find the box on
the outer boundary with the biggest label, breaking ties by preferring
the southernmost such box.  Let $a_{\binom{n}{2}}$ be the label of
this box, and $\alpha_{\binom{n}{2}}$ the row this box is in.  Now
remove the box and do jeu de taquin to fill the space, eventually
leaving a ``hole'' in the northwest corner of the tableau.  Repeat to
recover $a_{\binom{n}{2}-1}$ and $\alpha_{\binom{n}{2}-1}$ and so on
until the tableau is empty.  To adjust this so that it works for $w_n$
rather than the long word, we simply increase each $\alpha_k$ by 1.

\begin{prop}
The Edelman-Greene correspondence sends an rc-graph for $w_n$ to a
$Q$-tableau in which the label $i$ occurs only in rows $i-1$ and $i$.
\end{prop}

\begin{proof}
Since $Q$ is column-strict, the label $i$ cannot occur in any row
strictly below than row $i$.

Now suppose the label $i$ occurs in row $j$ for some $j<i-1$.  Then
the rightmost entry in row $j$ must have label $k$ for some $k\geq i$.
This entry will never be moved by the jeu de taquin during the
evacuation procedure, so, eventually, we will have an element
$(k,j+1)$.  But, since $j+1<k$, the $j+1$-st anti-diagonal does not
meet the $k$-th row, so we could not have started with an rc-graph
for $w_n$ in the first place.
\end{proof}

Now we can define a bijection from rc-graphs for $w_n$ to partitions
fitting inside $\delta_n$ by letting the partition associated to an
rc-graph be the boxes whose label matches the row number in the
$Q$-tableau corresponding to the rc-graph.  This is an injection
since the Edelman-Greene correspondence is injective, and this is
therefore a bijection due to proposition \ref{prop:count}.

\begin{prop}
The bijection given in Section \ref{sect:elem-bij} is the same as the
one given by the Edelman-Greene correspondence.
\end{prop}

\begin{proof}
For an rc-graph with $i$ cross pieces in the $k$-th row, the
Edelman-Greene corrspondence produces a tableau with $i$ occurrences
of the letter $k$, or, equivalently, a partition whose conjugate has
exactly $n-k+1-i$ parts of size $k-1$.  An rc-graph for any
permutation of $S_{n+1}$ has $n+1-k$ places, and therefore $n+1-k-i$
elbow joints, in the $k$-th row, so the two bijections are the same.
\end{proof}

Since the proof of proposition \ref{prop:elem-bij} gives a surjection
and the Edelman-Greene correspondence is known to be an injection, a
purely bijective proof omitting the counting lemma is possible.

\section{The Transposition Involution}\label{sect:tree-inv}

Since $w_n$ is its own inverse in $S_n$, transposing an rc-graph for
$w_n$ gives another rc-graph for $w_n$.  It is a natural question to
ask what this involution translates to on partitions.  Unfortunately,
the description of this involution on partitions is not immediately
evident.  However, described on bracketings of a string of length
$n+1$ subject to a binary nonassociative operation, it turns out to be
simply reversing the string along with the brackets.  Equivalently,
under the obvious bijection to binary trees, this corresponds to
flipping the tree along its vertical axis.

We describe the map from rc-graphs for $w_n$ to bracketings as
follows.  For convenience, let the ``letters'' of the string be the
numbers from $1$ to $n+1$.  Now, let $D$ be an rc-graph for $w_n$,
and for each $(i,j)\in E(D)$ (the set of locations of (nontrivial)
elbow joints in $D$), place a left bracket before the letter $j$ and a
right bracket after the letter $n+2-i$.  It is clear such a map is
injective and sends the transposition involution on rc-graphs to
reversal of order on parenthesizations.  What remains to be shown is
that this actually gives a proper full bracketing for a binary
associative operation.  Actually, more than this is true; the pair of
brackets aassociated with each elbow joint is in fact a matching pair.

We prove this by induction on the generalized inverse chute moves in
the proof of proposition \ref{prop:elem-bij}.  $D_\mathrm{bot}$
corresponds to the bracketing $(1(2(\cdots(n n+1)\cdots)))$, and each
elbow joint clearly corresponds to a matching pair of brackets.  Now
suppose there is valid generalized inverse chute move moving an elbow
joint at $(i,j)$ to $(i^\prime,j^\prime)$.  The elbow joint at $(i,j)$
corresponds to a matching pair of brackets with the left bracket
before the letter $j$ and the right bracket after the letter $n+2-i$.
The second condition for a valid generalized inverse chute move forces
the next right bracket to also occur between the letters $n+2-i$ and
$n+3-i$; the first condition forces us to have another left bracket to
the left of the letter $j$ matching a right bracket after the letter
$n+2-i^\prime$, although this could be an imaginary pair of brackets
around the letter $j$ (corresponding to a trivial required elbow joint
at $(j,n+2-j)$).  The remaining conditions merely state that the
original pair of brackets is a matching pair, which for us is true by
induction.  We can draw the situation as follows:

$$(1\cdots(j^\prime ((j \cdots n+2-i^\prime) \cdots n+2-i))\cdots n+1)$$

The generalized inverse chute move shifts the parentheses to the
following configuration:

$$(1\cdots((j^\prime (j \cdots n+2-i^\prime)) \cdots n+2-i)\cdots n+1)$$

Clearly, the new pair of brackets is a matching pair whenever we
start with a proper full bracketing.  Translated into the language of
binary trees, this operation is (left) rotation, an operation used in
many schemes for keeping binary search trees balanced.

\section{Multiplicity on $X_{w^\prime_n}$}\label{sect:geom}

A {\it (complete) flag} $\mathcal{F}$ in $\mathbb{C}^n$ is a sequence
of subspaces $\{0\}=F_0\subset F_1\subset F_2\subset\cdots\subset
F_n=\mathbb{C}^n$ such that the subspace $F_i$ has dimension $i$.
Fixing a basis for $\mathbb{C}^n$, we can represent $\mathcal{F}$
non-uniquely by an invertible matrix $m\in GL(n)$, where the first $i$
columns of $m$ form a basis for $F_i$.  Two matrices $m$ and
$m^\prime$ represent the same flag precisely when $m^\prime=m\cdot b$
for some $b\in B$, the group of upper-triangular matrices; as a result
the {\em flag variety} which parameterizes set of all flags is the
quotient $GL(n)/B$.  Note that $GL(n)$ and its subgroups $B$ and
$B_-$, the lower triangular matrices, act on $GL(n)/B$ on the left.
Given our choice of basis, $GL(n)/B$ has a distinguished flag $E$
called the {\em standard flag} represented by the identity matrix.

For each permutation $w\in S_n$ there is a subvariety $X_w$ of
$GL(n)/B$ known as the {\em Schubert variety}; $X_w$ is the closure of
the left $B$-orbit of the flag $F=w(E)$.  $F$ is the flag whose $i$-th
vector space $F_i$ is spanned by the vectors $e_w(1),\ldots,e_w(i)$,
or, alternatively, $F$ is the flag represented by $w$ as a permutation
matrix.  Given two permutations $v$ and $w$, $X_v \subseteq X_w$ iff
$v<w$ in the Bruhat order on $S_n$; in particular, for $e$ the
identity permutation, $X_e\in X_w$ for every $w\in S_n$.  Note that
$X_e$ consists of a single point, namely the flag $E$.  A dense open
neighborhood of $E$ is given by $\Omega^\circ_e$, the orbit of $E$
under the left action of $B_-$.  A detailed treatment of flag and
Schubert varieties can be found in, for example, \cite{Ful97}.

The {\it multiplicity} of a variety $X$ at a point $p$ is the degree
of the projective tangent cone
$\operatorname{Proj}(\operatorname*{gr}_{\mathfrak{m}_p}
\mathcal{O}_{X,p})$, considered as a subvariety of the projective
tangent space $\operatorname{Proj} (\operatorname*{Sym}^*
\mathfrak{m}_p/\mathfrak{m}^2_p)$.  The multiplicity is one measure of
``how singular'' $X$ is at $p$; in particular it is always 1 if $X$ is
smooth at $p$.  General semi-continuity theorems imply that the
multiplicity of $X_w$ at $X_e$ is at least the multiplicity of $X_w$
at any other point.

The multiplicity of $X_w$ at $X_e$ can be calculated using local
equations for $X_w$ on $\Omega^\circ_e$.  In general, these equations
will be a specialization of the equations for the matrix Schubert
varieties given in \cite{Ful92}.  However, if $w$ satisfies the
condition that, for every $(i,j)$ with $i+j>n$, either
$(w_0w)^{-1}(i)\leq j$ or $w_0w(j)\leq i$, the local equations are
exactly the equations for the matrix Schubert varieties.  Therefore,
in this special case, the multiplicity of $X_w$ at $X_e$ is exactly
the degree of the matrix Schubert variety.

Let $w^\prime_n$ be the permutation $(n+2)23\cdots(n+1)1\in S_{n+2}$.
These permutations satisfy the condition stated above, so the
multiplicity of $X_{w^\prime_n}$ at $X_e$ is given by the degree of
the matrix Schubert variety.  Knutson and Miller \cite{KM04} relate
rc-graphs to degenerations of matrix Schubert varieties so that, in
particular, the degree of a matrix Schubert variety is given by
$\mathfrak{S}_{w_0w}(1,\ldots,1)$.

Note that $w_0w^\prime_n=1(n+1)\cdots2(n+2)$, and Schubert polynomials
are unchanged under the inclusion of $S_k$ into $S_k+1$ fixing the
last element, as can easily be seen by adding an anti-diagonal of
elbow joints to every rc-graph.  Therefore, the multiplicity of
$X_{w^\prime_n}$ at $X_e$ is given by
$\mathfrak{S}_{w_n}(1,\cdots,1)=C_n$.

\section{Acknowledgements}

Thanks to Ezra Miller for suggesting this project and getting me
started, and Bernd Sturmfels, Brian Rothbach, and Mark Haiman for
useful suggestions along the way.  Thanks also to Alexander Yong for
convincing me to finally write up these results, and again to Ezra
Miller for the LaTeX macros for drawing rc-graphs.

\end{document}